\documentclass{article}
\usepackage{fullpage}
\def\C{{\rm C \kern-.48em\vrule width.06em height.57em depth-.02em \kern.48em}}
\def\Z{{{\rm Z}\kern-.28em{\rm Z}}}
\def\N{{{\rm I}\kern-.16em{\rm N}}}
\def\R{{{\rm I}\kern-.16em{\rm R}}}
\def\Re{\mathop{\rm Re}\nolimits}
\def\Im{\mathop{\rm Im}\nolimits}
\def\deg{\mathop{\rm deg}\nolimits}

\def\barz{\overline{z}}

\def\ii{{\rm i}}
\def\ee{{\rm e}}
\def\fromto{\mathop{{..}}}
\def\eqbd{\mathop{{:}{=}}}
\def\bdeq{\mathop{{=}{:}}}
\newtheorem{thm}{Theorem}
\newenvironment{@abssec}[1]{%
     \if@twocolumn
       \section*{#1}%
     \else
       \vspace{.05in}\footnotesize
       \parindent .2in
         {\bfseries #1. }\ignorespaces
     \fi}
     {\if@twocolumn\else\par\vspace{.1in}\fi}
\newcommand {\proof} {\par{\bf Proof}. \hskip 0.4cm}
\newcommand {\eproof}
      {\hfill
        {\ \vbox{\hrule\hbox{\vrule height1.3ex\hskip0.8ex\vrule}\hrule}}
        \vskip 0.4cm \par}
\newenvironment{keywords}{\begin{@abssec}{Key words}}{\end{@abssec}}
\newenvironment{AMS}{\begin{@abssec}{AMS subject classification}}{\end{@abssec}}
\begin{document}

\title{Hermite-Biehler, Routh-Hurwitz, and total positivity}
\author{Olga Holtz\footnote{On leave from Department of Computer Sciences,
 University of Wisconsin, Madison, Wisconsin 53706 U.S.A. 
Supported by Alexander von Humboldt Foundation.}  \\ \\
         Institut f\"ur Mathematik, MA 4-5 \\
         Technische Universit\"at Berlin \\
         D-10623 Berlin, Germany \\
         holtz@math.tu-berlin.de } 
\date{}
\maketitle


\begin{abstract} \noindent
Simple proofs of the Hermite-Biehler and Routh-Hurwitz theorems
are presented. The total non-negativity of the Hurwitz matrix of a stable 
real polynomial follows as an immediate corollary.
\end{abstract} 

\section{Introduction}

The classical result of Routh-Hurwitz on the stability of polynomials
is now more than a century old. Its rich connections with other areas 
of analysis and algebra have been exposed in many subsequent works.
Monographs by Postnikov~\cite{P} and by Chebotarev and 
Meiman~\cite{CM} give a detailed account of such related questions, 
including the amplitude-phase interpretation of stability, Sturm
chains, Cauchy indices, the principle of the argument, continued fractions, 
Hermite-Biehler theorem, and rational lossless functions. The interested 
reader should also see Barnett and Siljak's centennial survey~\cite{BS} 
and references therein to find out what control theory problems can be 
solved using Routh's algorithm. Krein and Naimark's survey~\cite{KN} and 
Householder's article~\cite{H}  explore connections 
of the Routh-Hurwitz scheme with Bezoutians, while Genin's article~\cite{G} 
emphasizes connections with Euclid's algorithm and orthogonal polynomials,
and presents a generalized Routh-Hurwitz algorithm suitable, e.g.,
for testing nonnegativity of a polynomial.  Asner~\cite{AS} and 
Kemperman~\cite{KE} found the link between stability and total 
positivity of the Hurwitz matrix. The Routh-Hurwitz algorithm,
originally formulated for real polynomials, has been extended to
complex polynomials (see~\cite{GA} or~\cite{F}) and further to wider 
classes of analytic functions (see~\cite{CM}). This list is in no way 
exhaustive, since the existing literature on the subject is enormous.

This note is not a survey of the field. Nor does it present a new approach 
to the subject. The note serves to derive, in a most elementary and 
economical way, three basic results in the Routh-Hurwitz theory, namely, 
the Hermite-Biehler theorem, the Routh-Hurwitz criterion, and the total 
positivity of the Hurwitz matrix. Because of this last issue, 
the setup is restricted to real polynomials. However, the proof 
of the Hermite-Biehler theorem extends verbatim to the complex case 
by considering the (generalized) odd and even part of a polynomial.
My approach is minimalistic. For example, orthogonality of polynomials, 
rational lossless functions and the like are not discussed when the only 
fact needed is interlacing of roots. The point is not that these connections 
are unimportant, but that they are not needed for a quick and direct 
derivation of the basics of the Routh-Hurwitz theory. 

Two papers,~\cite{AD} and~\cite{Me}, offer alternative elementary proofs of 
the Routh-Hurwitz scheme (labeled Theorem~\ref{rh} in this note). The proof 
in~\cite{AD} is based on geometric considerations in the complex plane and 
the proof in~\cite{Me}, simpler in my opinion, solely on continuity of the 
roots of a  polynomial.  Here, on the other hand, three results are derived
essentially at once, Routh-Hurwitz being a direct consequence of root
interlacing obtained in Hermite-Biehler, and the total nonnegativity of
the Hurwitz matrix a direct consequence of Routh-Hurwitz interpreted as
a matrix factorization formula.

\section{Proofs}

{\bf Definition}. \hskip 0.4cm
 A polynomial $f$ is {\em stable\/} if the condition $f(z)=0$ implies $\Re z<0$.

The following is a version of the Hermite-Biehler theorem (\cite{HM},~\cite{BL}).

\begin{thm} \label{hb} Let $f(x)=p(x^2)+xq(x^2)$, $f(x)$, 
$p(x)$, $q(x)\in \R[x]$. The following are equivalent.
\begin{description}
\item{\bf A}. The polynomial $f$ is stable.  

\item{\bf B}. The polynomials $p(-x^2)$ and $xq(-x^2)$ have simple real interlacing 
      roots and $\Re {p(z_0^2) \over z_0 q(z_0^2)} >0$ for some $z_0$ with $\Re z_0>0$.
\end{description}
\end{thm}

\proof A$\Longrightarrow$B:  If $f$ is stable, then  $f(z) = a\prod_j(z-z_j)$ with all
$z_j$ in the left half-plane. If $\Im z>0$, then  $|\ii z + \barz_j| > |\ii z - z_j|$ for
all $j$, hence $|f(\ii \overline{z})|>|f(\ii z)|$ or, by expanding the squares of both
absolute values and simplifying, $\Im p(-z^2)\overline{z} q(-\overline{z}^2)<0$. 
This implies that the functions 
\begin{equation} 
z\mapsto  {p(-z^2) \over z q(-z^2)}, \qquad z\mapsto  {z q(-z^2) \over p(-z^2)} \label{func} 
\end{equation}
take on real values only on the real axis. Hence any non-trivial real linear combination 
\begin{equation} \lambda p(-z^2)+\mu z q(-z^2), \qquad \lambda^2+\mu^2 \neq 0, \label{comb}
 \end{equation}
has only real roots. Next, $\gcd( p(-x^2),xq(-x^2))=1$, for if not, then $f$ would have either 
two roots with opposite real parts or one on the imaginary axis.  So, if~(\ref{comb}) had a 
multiple root, one of the functions~(\ref{func}) would have a high-order crossing with some 
horizontal line.  But if $g(x)-r=(x-x_0)^k h(x)$, $h(x_0)\neq 0$, for analytic functions $g$, 
$h$, and $k>1$, then, for small $\varepsilon>0$, the equation $g(x)=r-\varepsilon^k$ has solutions
 $x_0+\ee^{\ii\pi(1+2j)/k} h(x_0)^{-1/k} \varepsilon+ o(\varepsilon)$, $j=1,\ldots, k$. Hence the 
function $g$ takes on real values somewhere off the real axis.
This shows that no combination~(\ref{comb}) has a multiple root.  This also implies that
the roots of $p(-x^2)$ and $xq(-x^2)$ interlace, for if not, then one of the functions would 
preserve its sign on the interval between two consecutive roots of the other, hence, by a 
standard argument, there would be a  combination~(\ref{comb}) with a multiple root inside 
that interval.

B$\Longrightarrow$A:  If Condition B holds, then the function 
$ z \mapsto \Re \left( {p(z^2) \over z q(z^2)} \right)$ does not change its sign 
in the half-plane $\Re z>0$ and that sign is positive. So, the equation ${p(z^2) \over z q(z^2)} +1=0$
or, equivalently, $f(z)=0$, has no solution with $\Re z>0$. The roots of $p(-x^2)$ and $xq(-x^2)$ 
are distinct, so there is no solution to $f(z)=0$ on the imaginary axis either.\eproof
  
{\bf Remark}.\hskip 0.4cm The beginning of this proof is in the spirit of the argument 
from~\cite[pp.~13--15]{CM}. demonstrating that A implies that 
$$ \Im \left( {p(-z^2) \over z q(-z^2)}\right) 
      \Im z <0 \qquad {\rm whenever} \quad  z \notin \R. $$

The following Theorem is the essence of the Routh-Hurwitz scheme. It is proved in monographs 
using  Cauchy indices, Sturm  chains or the principle of the argument (see, 
e.g.,~\cite[pp.~225--230]{GA}). A nice elementary proof is given in~\cite{Me}.
Here is a different elementary argument based on Theorem~\ref{hb}.

\begin{thm} \label{rh} The polynomial $f(x)=p(x^2)+xq(x^2)$, with 
$p(x)$, $q(x)\in \R[x]$, is stable if and
only if $c\eqbd p(0)/q(0)>0$ and the polynomial $\widetilde{f}(x)=
\widetilde{p}(x^2)+x\widetilde{q}(x^2)$ is stable, where 
$\widetilde{p}(x) \eqbd q(x)$, $\widetilde{q}(x) \eqbd {1\over x} 
(p(x)-cq(x))$.
\end{thm}

\proof Necessity. Condition B in Theorem~\ref{hb} is 
equivalent to $p$ and $q$ satisfying $p(0)q(0)>0$ and having only 
simple zeros, all negative,  interlacing, the rightmost zero
being that of $p$.

Let the pair $(p,q)$ satisfy Condition B, let $x_n<\cdots<x_1$ be the zeros of $p$, 
and $y_k<\cdots<y_1$ the zeros of $q$, and assume wlog that $p(0)>0$.
Then, with $y_n$ any point to the left of $x_n$ in case $k=n-1$, we have $p$
and $q$ of opposite sign in $(y_j\fromto x_j)$, all $j$, and also
$(-1)^jp(y_j)>0$ for all $j$. But then, for any $c\ge0$, the polynomial
$r:= p - c q$ has the same sign as $p$ on $[y_j\fromto x_j)$, all $j$. In
particular, also $(-1)^jr(y_j)>0$, all $j$, and this implies that, in each
of the $n-1$ intervals $(y_{j+1}\fromto y_j)$, $r$ has an odd zero. If now,
specifically, $c = p(0)/q(0)$ (which is positive, by assumption), then $r$ 
also
has a zero at 0, and since its degree is no bigger than $n$, those $n-1$ odd
zeros must all be simple. But this implies that
$\widetilde q $ is of degree $n-1$, with all its zeros simple and
negative, and these zeros separate those of $\widetilde p := q$, and, in
particular, $q(y_1)$ has the sign opposite to $r(y_1)$, i.e., to $p(y_1)$,
i.e., is positive, hence both $\widetilde p$ and $\widetilde q$ are positive at $0$.
In short, if $(p,q)$ satisfies Condition B of Theorem~\ref{hb}, then so does
the pair $(\widetilde p,\widetilde q)$.

Sufficiency. Suppose $\widetilde{f}(x)$ is stable and $c>0$. 
Since $p(x^2)=c\widetilde{p}(x^2)+x^2\widetilde{q}(x^2)$, $q(x^2)=\widetilde{p}(x^2)$, and,
by  Theorem~\ref{hb}, 
$\Re \left( {\widetilde{p}(z^2) \over z \widetilde{q}(z^2)} \right)>0 $ 
whenever $\Re z>0$, one obtains
$$\Re \left( {p(z^2)\over zq(z^2)} \right) =
\Re \left( {c \over z} \right)+ \Re \left( {z \widetilde{q}(z^2) \over 
\widetilde{p}(z^2)} \right)> \Re \left( {c\over z} \right)>0 \qquad
{\rm whenever} \quad \Re z>0.$$ 
 Finally, if $\widetilde{p}(x^2)$ and $x \widetilde{q}(x^2)$ 
are relatively prime, so are $p(x^2)$ and $x q(x^2)$. This proves
that Condition B of Theorem~\ref{hb} is met. \eproof

Theorem~\ref{rh} implies the following version of the 
Routh-Hurwitz theorem (\cite{RT}, \cite{HW}).

\begin{thm} \label{rhnew} The polynomial  $f(x)\bdeq \sum_{j=0}^n {a_j} x^j$ with
$a_0>0$ is stable if and only if its infinite Hurwitz matrix $H(f)$ is a 
product of the form 
\begin{equation} 
H(f)=J(c_1)\cdots J(c_n) H(b), \label{fctr}
\end{equation}
 with all parameters $c_j$, $j=1, \ldots, n$, positive, and
$b$ a positive polynomial of degree 0.  Here  
$$ H(f) \eqbd \left( \begin{array}{ccccc} a_0 & a_2 & a_4 & a_6 & \cdots \\
     0 & a_1 & a_3 & a_5 & \cdots \\
     0 & a_0 & a_2 & a_4 & \cdots \\
     0 &  0 & a_1 & a_3 &  \cdots \\
    \vdots & \vdots & \vdots & \vdots & \ddots \end{array} \right),
\quad J(c)\eqbd \left( \begin{array}{cccccc} c & 1 & 0 & 0 & 0 & \cdots \\
     0 & 0 & 1 & 0 & 0 & \cdots \\
     0 & 0  & c & 1 & 0 & \cdots \\
     0 &  0 & 0 & 0 & 1 & \cdots \\
     0 &  0 & 0 & 0 & c & \cdots \\ 
    \vdots & \vdots & \vdots &  \vdots & 
    \vdots & \ddots \end{array} \right). $$
\end{thm}

\proof Prove by induction that a polynomial $f$ of degree $n$, $f(0)>0$,
 is stable if  and only if the first $n+1$ leading principal minors 
$\Delta_j(f)$, $j=1, 
\ldots, n+1$, of $H(f)$ are positive and the factorization~(\ref{fctr})
holds. Indeed, if $n=0$, both properties are valid trivially. If
$n>0$, then, by the Lemma, $f$ is stable if and only if $c>0$ and 
$\widetilde{f}$
is stable. But, as one readily verifies, $H(f)=J(c)H(\widetilde{f})$, hence, 
in  particular, 
$\Delta_{j+1}(f)=c\widetilde{f}(0)\Delta_j(\widetilde{f})$, $j=0, 1, \ldots$;
here $\Delta_0 \eqbd 1$.  Since $\deg \widetilde{f}=\deg f -1$, 
$\widetilde{f}$  satisfies the inductive hypothesis, hence so does $f$. 
\eproof

\begin{thm} \label{tp} The Hurwitz matrix of a stable polynomial $f$ satisfying
$f(0)>0$ is totally nonnegative. 
\end{thm}

\proof By Theorem~\ref{rhnew}, the factorization~(\ref{fctr})
holds with all parameters positive. By inspection, each factor is totally 
nonnegative, hence their product $H(f)$ is also totally nonnegative.
\eproof

Theorem~\ref{tp} was first proved in~\cite{AS} and~\cite{KE}.

\section*{Acknowledgements} I am grateful to Carl de Boor and Hans Schneider for 
their detailed suggestions for the improvement of the manuscript and to the referee 
for helpful critique.


\begin{thebibliography}{ll}
\bibitem{AD}
J. J. Anagnost and C. A. Desoer, {\em 
An elementary proof of the Routh-Hurwitz stability criterion, \/}
Circuits Systems Signal Process. 10 (1991), no. 1, 101--114.

\bibitem{AS} B. A. Asner, {\em On the total nonnegativity of the Hurwitz matrix,\/}
SIAM J. Appl. Math., 1970, 18, 407--414.

\bibitem{BL} M. Biehler, {\em Sur une classe d'\'{e}quations alg\'ebriques  dont toutes les racines sont r\'eelles. \/} J. reine angew. Math., 1879, 87, 350--352.

\bibitem{BS} S. Barnett, D. D. Siljak, {\em Routh's algorithm: a
centennial survey,\/} SIAM Review, 1977, 19(3), 472--489.

\bibitem{CM} N. G. Chebotarev and N. N. Meiman. 
The Routh-Hurwitz problem for polynomials and entire functions. Real quasipolynomials
with $r=3$, $s=1$. (Russian) Appendix by G. S. Barhin and A. N. Hovanskii. 
Trudy Mat. Inst. Steklov. 26, (1949). 

\bibitem{F} E. Frank, {\em On the zeros of polynomials with complex 
coefficients,\/} Bull. Amer. Math. Soc., 1946, 52, 144--157.

\bibitem{GA} F. R. Gantmacher.  The Theory of Matrices, vol.II; New York,
Chelsea Publ., 1959.

\bibitem{G} Y. V. Genin, {\em Euclid algorithm, orthogonal polynomials, and
generalized Routh-Hurwitz algorithm,\/} Linear Algebra Appl., 1996, 246: 
131--158. 

\bibitem{HM} C. Hermite, {\em Sur les nombre des racines d'une \'{e}quation alg\'{e}brique
comprise entre des limites donn\'{e}es,\/} J. reine angew. Math., 1856, 52, 39--51.

\bibitem{H} A. S. Householder, {\em Bezoutians, elimination and localization,
\/} SIAM Review, 1970, 12 (1), 73--78. 

\bibitem{HW} A. Hurwitz, {\em  \"{U}ber die Bedingungen, unter welchen eine Gleichung 
nur Wurzeln mit negativen reellen Teilen besitzt,\/} Math. Ann., 1895, 46, 273--284.

\bibitem{KE} J. H. B. Kemperman, {\em A Hurwitz matrix is totally positive,\/} SIAM J. Math.
Anal., 1982, 13, 331--341.

\bibitem{KN} M. G. Krein and M. A, Naimark, {\em The method of symmetric
and Hermitian forms in the theory of the separation of the roots of 
algebraic equations,\/} Linear and Multilinear Algebra, 1981, 10, 265--308.

\bibitem{Me}
G. Meinsma, {\em  Elementary proof of the Routh-Hurwitz test, \/}
Systems Control Lett. 25 (1995), no. 4, 237--242.

\bibitem{P} M. M. Postnikov. Stable polynomials. (Russian) Moscow,
Nauka, 1981.

\bibitem{RT} E. J. Routh. Stability of a given state of motion. London, 1877.

\end{thebibliography}
\end{document}